

\baselineskip=14pt
\parskip=10pt

\font\eightrm=cmr8 

\magnification=\magstephalf

\def\1{{\overline{1}}}
\def\2{{\overline{2}}}
\parindent=0pt
\overfullrule=0in

\def\frac#1#2{{#1 \over #2}}
\centerline
{\bf  
Counting Standard Young Tableaux  With Restricted Runs
}
\bigskip
\centerline
{\it Manuel KAUERS and Doron ZEILBERGER}
\bigskip
\quad\quad\quad\quad {\it Dedicated to super-enumerators Ian Goulden and David Jackson\footnote{$^1$}
{\eightrm This article was originally submitted to the electronic journal ``Algebraic Combinatorics''
founded by Goulden and Jackson, following a solicitation  for a special issue
in honor of Goulden and Jackson. On July 27, 2020 we got an email message from one of the
editors-in-chief, Akihiro Munemasa, informing us that, after an initial review, it is ``unlikely
to meet the standards of depth and originality that the journal is seeking''. Consequently this
article will remain in the `Personal journal of Shalosh B. Ekhad and Doron Zeilberger', the homepage of Manuel Kauers, and
of course arxiv.org. Let the readers decide about its depth and originality.}}

{\bf Abstract.} The number of standard Young tableaux whose shape is a $k$ by $n$ rectangle is
famously $\frac{(nk)!\,0!\,1!\, \cdots (k-1)!}{(n+k-1)!(n+k-2)! \cdots n!}$, implying
that for each specific $k$, that sequence satisfies a linear recurrence equation with polynomial
coefficients of the {\bf first} order. But what about counting standard Young tableaux where
certain ``run lengths" are forbidden? Then things seem to get much more complicated.
In this tribute to the legendary enumerative pair {\it Goulden \& Jackson} we
investigate these intriguing sequences, and conjecture that if the number of rows is larger than two,
then these sequences are generally {\bf not} $P$-recursive. On the positive side, we conjecture
that these sequences have `nice' asymptotic behavior.
We pledge donations to the OEIS in honor of the first solvers of these conjectures.

{\bf Preface}

Some combinatorial families are easy to count, for example the number of subsests of an $n$-element set,
that can be computed in {\it logarithmic} time (in base $2$). Also easy is counting the number of permutations,
that can be computed in linear-time. Then you have the really hard ones, for example the number
of $n \times n$ Latin squares  and the number of self-avoiding walks of length $n$, for which
we will probably never know the exact value of the $1000$-th term.

Both the number of subsets of an $n$-element set, $2^n$, and the number of permutations, $n!$,
as well as the famous Catalan numbers $(2n)!/(n!(n+1)!)$ (OEIS sequence $A108$) are said to have a {\bf closed-form} formula.
They satisfy  a {\it first-order} linear recurrence equation with polynomial coefficients
$$
a(n+1)-2a(n)=0 \quad , \quad a(n+1)-na(n)=0 \quad, \quad  (n+2) a(n+1) - 2 (2n+1) a(n)=0 \quad .
$$
(Note that the first equation is even better, it is {\it constant coefficients}).

Many natural families satisfy the next-best thing to being closed-form, they satisfy 
a {\bf linear-recurrence equation} with {\bf polynomial} coefficients, but not necessarily
of first order. Such sequences, called $P$-recursive, or {\it holonomic} (see [KP]), satisfy
an equation of the form
$$
\sum_{i=0}^{L} p_i(n) a(n+i) \, = \, 0 \quad ,
$$
for some positive integer $L$ and some {\bf polynomials} in $n$, $p_0(n), \dots , p_L(n)$.

The most famous such sequences that are {\bf not} closed-forms are the
Fibonacci numbers, $F_n$ (OEIS sequence $A45$), the number of involutions of an $n$-element set
(permutations that are equal to their inverse) $w_n$ (OEIS sequence $A85$), and the Motzkin numbers, $M_n$,
the number of words of length $n$ in the alphabet $\{0,-1,1\}$ that add-up to $0$ and all whose
partial sums are non-negative (OEIS sequence $A1006)$. They satisfy, respectively, the recurrences 
$$
a(n+2)-a(n+1)-a(n)=0 \quad , \quad  a(n+2)- a(n+1) -(n+1)a(n)=0 \quad, \quad
$$
$$
(n+4)a(n+2) - (2n+5)a(n+1) - (3n+3)a(n)=0 \quad .
$$

Given a natural combinatorial family, parameterized by $n$, it is very interesting to know whether or not
the sequence of integers that enumerates it happens to
be $P$-recursive. This is interesting both conceptually and computationally, since
a linear recurrence makes it easy to compute many terms, as well as deriving the
asymptotics. Sometimes proving that a given sequence, or family of sequences, is (are) $P$-recursive
is highly non-trivial, see for example [GJR], [GJ], and [Ge2].

In this modest tribute to Ian Goulden and David Jackson we will raise the question whether
a certain very {\it natural\/} family of combinatorial sequences is $P$-recursive, and
give ample computational evidence that generally they are, {\it probably}, {\bf not}.

{\bf Maple packages} 

This article is accompanied by the Maple packages {\tt YoungT.txt} and {\tt Tableaux3R.txt}, available from 
{\tt http://www.math.rutgers.edu/\~{}zeilberg/mamarim/mamarimhtml/cyt.html} \quad .

{\bf Counting Standard Young Tableaux with Restricted Runs}

Recall that a partition of a positive integer $m$ is a weakly-decreasing list of positive integers
$\lambda=[\lambda_1, \dots, \lambda_k]$, that add-up to $m$, also called a {\it shape}, and a
{\it standard Young tableau} of shape $\lambda$  is a left-justified array of $k$ rows
with $\lambda_i$ boxes in the $i$-th row, where the integers $\{1, \dots , m\}$
are filled in the boxes so that both rows and columns are  increasing.
For example, here is a standard Young tableau of shape $[3,3,2]$
$$
\matrix{ 1 & 3 & 4 \cr
         2 & 5 & 7 \cr
         6 & 8} \quad,
$$

and here is one of shape  $[5,5,4]$

$$
\matrix{ 1 & 3 & 4 & 6 & 7\cr
         2 & 5 & 8 &9  &11\cr
         10 & 12 & 13 &14} \quad .
$$

The number of standard Young tableaux of shape $[\lambda_1, \dots, \lambda_k]$ is famously given by the {\bf Young-Frobenius} formula
(equivalent to the {\it hook-length formula})
$$
\frac{(\lambda_1 + \dots + \lambda_k)!}{(\lambda_1+k-1)! \cdots \lambda_k!} \cdot \prod_{1 < i< j \leq k}  (\lambda_i-\lambda_j+j-i) \quad .
$$

In particular,  setting $\lambda_1= \dots = \lambda_k=n$, we get that
the number of standard Young tableaux whose shape is a $k$ by $n$ rectangle is given by
$$
\frac{(nk)!\,  0!\,1!\, \cdots (k-1)!}{(n+k-1)!(n+k-2)! \, \cdots \,n!} \quad .
$$
It follows that for each {\bf fixed} $k$, this sequence, in $n$, is $P$-recursive, in fact it even satisfies
a {\bf first-order} linear recurrence.

Let's define a {\bf run} in a standard Young tableau to be a {\bf maximal} string of {\bf consecutive integers}.
For example, in the following tableau of shape $[5,5,5]$

$$
\matrix{ 1 & 3 & 4 & 6 &   7\cr
         2 & 5 & 8 &9  &   11\cr
         10 & 12 & 13 & 14 & 15} \quad,
$$
we have

$\bullet$ First row: one run of length $1$,  ($1$), and two runs of length $2$ ($34$ and $67$)  ;

$\bullet$ Second row: three runs of length $1$, ($2$, $5$ and $11$)  and one run of length $2$ ($89$) ;

$\bullet$ Third row: one run of length $1$, ($10$)  and one run of length $4$ ($12,13,14,15$) .

We are interested in the following general question. Fix $k>1$.
Given {\it arbitrary} finite 
sets of positive integers $A_1, \dots, A_k$, (or infinite arithmetical progressions), compute the following integer sequence,
let's call it
$$
G_{A_1, \dots, A_k} (n) \quad,
$$
defined as the number of standard Young tableaux of shape $(n, \dots ,n)$ such that in each row $i$ , $1 \leq i \leq k$, 
none of the runs belongs to $A_i$. Of course if all the $A_i$'s are the empty set, we are back
to counting unrestricted standard Young tableaux, for which there is a nice closed-form formula, and
of course it is $P$-recursive.

The case of {\bf two} rows can be shown [EZ] to always give $P$-recursive sequences (in fact even something stronger is
true: the generating functions are algebraic formal power series). This is the case since $2\times n$ standard Young tableaux
are in easy bijection with {\it Dyck paths} of semi-length $n$.

More generally, as is well-known, and immediate to see,
standard Young tableaux are in easy bijection with {\bf lattice paths}. A $k$-rowed standard Young tableau
of shape $[\lambda_1, \dots, \lambda_k]$ corresponds to a lattice path in the $k$-dimensional hyper-cubical lattice,
from the origin $[0,\dots,0]$ to the point $[\lambda_1, \dots, \lambda_k]$, with {\bf unit positive steps}
${\bf e}_i:=[0^{i-1},1,0^{k-i}]$, that {\bf always} stay in the region $x_1 \geq x_2 \geq \dots \geq x_k$.
Given a standard Young tableau, the corresponding path is obtained by executing the step ${\bf e}_i$ at the $m$-th step,
if $m$ is located at the $i$-th row. So $G_{A_1, \dots, A_k}(n) $  is also the number
of $k$-dimensional lattice paths from the origin to $[n, \dots, n]$,
always staying in $x_1 \geq \dots x_k \geq 0$,  such that the walker
never has a run-length parallel to the $i$-th axis that belongs to the set $A_i$.

The same question makes sense for general walks, not necessarily those confined to $x_1 \geq x_2 \geq \dots \geq x_k$. 
It turns out that for this analogous question the sequences are {\bf always} $P$-recursive, as we will now
show.

{\bf Counting Lattice Walks with Restricted Runs}

In order to motivate the general case,
let's first give yet another proof, a bit more complicated than the usual one,
of the very easy fact that the generating function for the number of {\bf all} walks,
without restrictions, is given by the generating function
$$
\frac{1}{1-x_1 - \dots - x_k} \quad .
$$

Every walk corresponds to a {\bf word} in the alphabet $\{1, \dots k\}$, indicating which ${\bf e}_i$ it went through.
For example, the walk
$$
[0,0,0] \rightarrow [1,0,0] \rightarrow [1,0,1] \rightarrow [1,1,1] \rightarrow [1,2,1] \rightarrow [1,2,2] \rightarrow [1,3,2]
$$
corresponds to the word
$$
1 3 2 2 3 2 \quad.
$$
Given a word in $\{1,\dots, k\}$, we can write it in {\bf frequency notation} $b_1^{r_1} \dots b_l^{r_l}$,
where $b_{j+1} \neq b_{j}$, and $r_j \geq 1$.
for example, the above word $1 3 2 2 3 2$ is abbreviated $1^1 3^1 2^2 3^1 2^1$, and  the word
$1113322211$ is written $1^3 3^2 2^3 1^2$.

Let $F_i=F_i(x_1, \dots, x_k)$ be the {\bf weight-enumerator} of all words that end with the letter $i$.
Then obviously, for $i=1, \dots, k$
$$
F_i \, = \, \frac{x_i}{1-x_i} \left( 1+ \, \sum_{{{1 \leq j \leq k} \atop {j \neq i} }} F_j \right ) \quad.
$$

This is a system of $k$ linear equations with $k$ unknowns $F_1, \dots ,F_k$, whose solution is
easily seen to be given explicitly by
$$
F_i = \frac{x_i}{1-x_1 - \dots - x_k}   \quad .
$$
Finally the full generating function, $F$, is gotten by adding the weight of the {\bf empty} word, $1$, to the
sum of the $F_i$'s, getting 
$$
F= 1 + \sum_{i=1}^k F_i \quad,
$$
that implies the deep theorem
$$
F= \frac{1}{1-x_1 - \dots -  x_k} \quad .
$$

Note, in particular that the $F_i$, (and $F$) are {\bf rational functions} of the variables $x_1, \dots, x_k$.

To handle the restricted case, to find the weight-enumerator of all words in $1^{\lambda_1} \dots k^{\lambda_k}$ such that
when written in frequency notation $b_1^{r_1} \dots b_m^{r_m}$ we have that if $b_\alpha=i$ then $r_\alpha \not \in A_i$
(i.e. runs in the ${\bf e}_i$ direction can't be of a length that belongs to $A_i$), we have the modified system:
$$
F_i \, = \, \left ( \frac{x_i}{1-x_i} - \sum_{\beta \in A_i} {x_i}^{\beta} \right )  \left ( 1+ \, 
\sum_{{ {1 \leq j \leq k} \atop {j \neq i} } } F_j \right ) \quad.
$$
This is a system of $k$ linear equations  in the $k$ unknowns $F_1, \dots, F_k$,
with coefficients that are {\bf rational functions} in $x_1, \dots , x_k$.
Hence, by Cramer's rule, the $F_i$ are all rational functions of $x_1, \dots , x_k$,  and hence so is  $F=1+\sum_{i=1}^{k} F_i$.

Our sequence of interest is the sequence of coefficients of the {\bf  diagonal} of this rational function. Since
the diagonal of any formal power series that is a rational function is $D$-finite (see [Ge1][Z1][L][Z2]), it follows
that the sequence itself is $P$-recursive.

{\bf Back to Tableaux}

We strongly doubt that the multi-variable generating functions for restricted Young tableaux (alias restricted walks confined to $x_1 \geq \dots \geq x_k$)
are rational.
In order to explore these sequences, we need to generate as many terms as possible. Here is how to do it.
Let us fix $R_1, \dots, R_k$ and denote by
$g(\lambda_1, \dots, \lambda_k)$ the number of standard Young tableaux of shape $[\lambda_1,  \dots, \lambda_k]$
with no runs in the $i$-th row that belong to $R_i$, or
equivalently the number of walks in $x_1 \geq \dots \geq x_k \geq 0$,
from the origin to the point $[\lambda_1, \dots, \lambda_k]$ with no run-length in
the $x_i$-direction that belongs to $R_i$ (for $i=1, \dots, k$).  

In order to compute $g(\lambda_1, \dots, \lambda_k)$,
we need the more refined quantities ($1 \leq i \leq k$)
$g^{(i)}(\lambda_1, \dots, \lambda_k)$, that enumerate those walks that end with a step in the $x_i$-direction. 

We have the {\bf dynamic programming} recurrences ($1 \leq i \leq k$)
$$
g^{(i)}(\lambda_1, \dots, \lambda_k) \, = \,
\sum_{   {{1 \leq j \leq k} \atop {j\neq i} } }
\sum_{{{1 \leq r \leq \lambda_i} \atop {r \not \in A_i}} } g^{(j)}(\lambda_1, \dots, \lambda_{i-1}, \lambda_i-r, \lambda_{i+1}, \dots , \lambda_k) \quad,
$$
with the obvious initial conditions, and the {\bf boundary conditions}
$$
g^{(i)}(\lambda_1, \dots, \lambda_k) \, = \, 0  \quad,
$$
whenever $\lambda_1<\lambda_2$ or $\lambda_2<\lambda_3$, $\dots$, or $\lambda_{k-1}<\lambda_k$, or  $\lambda_k<0$.
Finally
$$
g(\lambda_1, \dots, \lambda_k)  = \sum_{i=1}^{k} g^{(i)}(\lambda_1, \dots, \lambda_k)  \quad .
$$
This is all implemented in the Maple package {\tt YoungT.txt} mentioned above.

\vfill\eject

{\bf Two Case Studies}

In spite of the fact that we were unable to think of a good reason why these sequences should be
$P$-recursive, we still hoped that they would be for a non-obvious reason. We focused on two
special cases to generate as many terms as we could.

$\bullet$ $G(n)$, the number of standard Young tableaux of shape $[n,n,n]$ where  each run, in each of the
three rows, must have length at least $2$. This is the case $A_1=A_2=A_3=\{1\}$ in the above notation.

$\bullet$ $H(n)$, the number of standard Young tableaux of shape $[n,n,n]$ where  all the run-lengths,
in each row are always odd. This is the case $A_1=A_2=A_3=\{2r+2: r \geq 0\}$ in the above notation.

Regarding $G(n)$, using the Maple package

{\tt http://www.math.rutgers.edu/\~{}zeilberg/tokhniot/Tableaux3R.txt}

we got that the sequence starts with (starting at $n=1$)
$$
0, 1, 1, 5, 15, 69, 304, 1518, 7807, 42314, 236621, 1364570, 8062975, 48680547, 299388670, 1871463427,  \dots  \quad .
$$
The first $200$ terms may be viewed here:

{\tt    https://sites.math.rutgers.edu/\~{}zeilberg/tokhniot/oTableaux3R1.txt} \quad .

The first $996$ terms are available here:

{\tt https://sites.math.rutgers.edu/\~{}zeilberg/tokhniot/CYT/GseqList.txt} \quad .

Regarding $H(n)$,
we got that the sequence starts with (starting at $n=1$)
$$
1, 2, 9, 46, 306, 2252, 18308, 158872, 1454570, 13888112, 137277741, 1396638636, 14561307281, 155040525128,  \dots \quad .
$$
The first $200$ terms can be viewed here:

{\tt    https://sites.math.rutgers.edu/\~{}zeilberg/tokhniot/oTableaux3R2.txt} \quad .

The first $965$ terms are available here: 

{\tt https://sites.math.rutgers.edu/\~{}zeilberg/tokhniot/CYT/HseqList.txt} \quad .

Even that many terms were not enough to guess a linear recurrence with polynomial coefficients, 
so if such a recurrence exists, it would be extremely complicated. But we can do better! The existence of a non-zero linear recurrence of
a given order and degree boils down to the existence of a non-zero solution to  a certain system of linear equations
with integer coefficients. If a non-trivial solution exists, then doing everything
modulo any prime would also have a solution. Conversely, if there is no solution modulo
that prime, there is no solution at all. Now we can generate many more terms,
and using the prime $p=45007$ we (or rather our computer) generated $5000$ terms, and
even these did not suffice. In other words if there exists such a recurrence of order $\leq K$ and
degree, in $n$, of the coefficients of degree $\leq K$, then $(K+1)^2+5 \geq 5000$, i.e. $K \geq 70$.

This leads us to make the following conjectures. One of us (DZ) is pledging a donation of $200$ US dollars
to the On-Line Encyclopedia of Integer Sequences (OEIS) in honor of the first prover, for each of the following four
conjectures.

{\bf Conjecture 1a}: The sequence $G(n)$  is  {\bf not} $P$-recursive.

{\bf Conjecture 1b}: The sequence $H(n)$  is  {\bf not} $P$-recursive.

Surprisingly, the asymptotics seems to be very nice. Using the nearly $1000$ terms in these sequences we are
safe in making the following conjectures.

{\bf Conjecture 2a}: There exists a constant $C_1$ (if possible, find it!) such that
$$
G(n) \, \asymp \, C_1 \frac{8^n}{n^4} \quad .
$$

We estimate $C_1$ to be close to  $0.521286$ .

{\bf Conjecture 2b}: There exists a constant $C_2$ (if possible, find it!) such that
$$
H(n) \, \asymp \, C_2 \frac{(7+ 5\, \sqrt{2})^n}{n^4} \quad .
$$

We estimate $C_2$ to be close to  $0.63892$.

This raises the more general question about these sequences. Is the asymptotics always of the form $C\mu^n\, n^\theta$ with
$\mu$ an algebraic number, and $\theta$ a rational number?

{\bf References}

[EZ] Shalosh B. Ekhad and Doron Zeilberger, {\it
Automatic counting of restricted Dyck paths via (numeric and symbolic) Dynamic Programming},
The Personal Journal of Shalosh B. Ekhad and Doron Zeilberger, June 3, 2020. \hfill\break
{\tt https://sites.math.rutgers.edu/\~{}zeilberg/mamarim/mamarimhtml/dyck.html}. \hfill\break
Also available from: {\tt https://arxiv.org/abs/2006.01961} \quad .

[Ge1] Ira Gessel, {\it Two theorems of rational power series},
Utilitas Math. {\bf 19 }(1981), 247-254.

[Ge2] Ira Gessel, {\it Counting Latin rectangles}, Bull. Amer. Math. Soc. (N.S.) {\bf 16} (1987), 79-82.

[GJR] Ian P. Goulden, David M. Jackson, and James W. Reilly,
{\it The Hammond series of a symmetric function and its application to P-recursiveness},
 SIAM J. Algebraic Discrete Methods {\bf 4} (1983), 179-193.

[GJ]  Ian P. Goulden and David M. Jackson,  {\it Labelled graphs with small vertex degrees and P-recursiveness},
 SIAM J. Algebraic Discrete Methods {\bf 7} (1986), 60-66.

[KP] Manuel Kauers  and Peter Paule, {\it ``The Concrete Tetrahedron''}, Springer, 2011.

[L] Leonard Lipshitz, {\it The diagonal of a D-finite power series is D-finite},
J. Algebra {\bf 113} (1988), 373-378.

[Z1] Doron Zeilberger, {\it Sister Celine's technique and its generalizations},
J. Math. Anal. Appl. {\bf 85} (1982), 114-145.

[Z2] Doron Zeilberger, {\it A holonomic systems approach to special functions identities},
J. Comput. Appl. Math. {\bf 32} (1990), 321-368.

\bigskip
\bigskip
\hrule
\bigskip
Manuel Kauers, Institute for Algebra, J. Kepler University Linz, Austria
E-mail {\tt manuel dot kauers at jku dot at}
\bigskip
Doron Zeilberger, Department of Mathematics, Rutgers University (New Brunswick), Hill Center-Busch Campus, 110 Frelinghuysen
Rd., Piscataway, NJ 08854-8019, USA. \hfill\break
Email: {\tt DoronZeil at gmail  dot com}   \quad .
\bigskip
\hrule
\bigskip
First Written: June 20, 2020.
This version: Aug. 8, 2020.
\bigskip
\hrule
\bigskip
Exclusively published in the Personal Journal of Shalosh B. Ekhad and Doron Zeilberger, Manuel Kauers' website, and arxiv.org .
\end